\documentclass[12pt]{amsart}
\usepackage{amscd,amssymb,graphics}
\theoremstyle{definition}

\theoremstyle{remark}

\numberwithin{equation}{section}

\def\norm#1{\left\Vert#1\right\Vert}
\def\R {{\mathbb R}}

\def\e{{\varepsilon}}
\def\Z {{\mathbb Z}}
\def\H{{\mathcal H}}

\def\ind{{\mathrm{ind}}}
\def\Ad{{\mathrm{Ad}}}

\begin{document}

\title[Questions of Eymard and Bekka]{On some
questions of Eymard and Bekka concerning amenability of homogeneous
spaces and induced representations}

\author[V. Pestov]{Vladimir Pestov}
\address{Department of Mathematics and Statistics, 
University of Ottawa, 585 King Edward Ave., Ottawa, Ontario, Canada K1N 6N5.
}
%\curraddr{}
\email{vpest283@science.uottawa.ca}
%\urladdr{http://www.mcs.vuw.ac.nz/$^\sim$vova}

\thanks{{\it 2000 Mathematical Subject Classification.} 
43A07, 43A65, 43A85, 22D30.}

\thanks{Partially supported by the Marsden Fund of the Royal Society of
New Zealand and by a University of Ottawa start-up grant.}

\keywords{}
\begin{abstract} 
Let $F\subseteq H\subseteq G$ be closed subgroups of a locally compact group.
In response to a 1972 question by Eymard, we construct an example
where the homogeneous factor space $G/F$ is amenable in the sense of 
Eymard--Greenleaf, while $H/F$ is not. (In our example, $G$ is discrete.) 
As a corollary which answers a 1990 question by Bekka, 
the induced representation
$\ind_H^G(\rho)$ can be amenable in the sense of Bekka even if $\rho$
is not amenable. The second example, answering another question by Bekka, 
shows that $\ind_H^G(\rho)$ need not
be amenable even if both the representation $\rho$ and the coset
space $G/H$ are amenable.
\\[3mm]
{\it R\'esum\'e.} Soient $F\subseteq H\subseteq G$ deux sous-groupes
ferm\'es d'un groupe localement compact $G$. 
En r\'eponse \`a une question pos\'ee
en 1972 par Eymard, nous construisons un groupe discret $G$
tel que l'espace homog\`ene
$G/F$ est moyennable au sens d'Eymard et Greenleaf, bien que 
$H/F$ n'est pas moyennable. On obtienne un corollaire qui r\'epond
\`a un probl\`eme pos\'e par Bekka en 1990: une 
repr\'esentation induite $\ind_H^G(\rho)$ peut \^etre moyennable
au sens de Bekka m\^eme si $\rho$ n'est pas moyennable.
Le deuxi\`eme exemple, qui r\'epond \`a une autre question de Bekka,
montre que $\ind_H^G(\rho)$ n'est pas n\'ecessairement moyennable m\^eme si
la repr\'esentation $\rho$ et l'espace homog\`ene $G/H$ sont 
l'un et l'autre moyennables.
\end{abstract}

\maketitle

\section{Introduction}
Let $H$ be a closed subgroup of a locally compact group $G$. The
homogeneous factor space $G/H$ is {\it amenable in the sense of Eymard and
Greenleaf} \cite{Ey,Gre}, if $L^\infty(G/H)$ supports a $G$-invariant mean.
If $H=\{e\}$, one obtains the
classical concept of an amenable locally compact group. 

A unitary representation $\rho$ of a group $G$ in a
Hilbert space $\H$ is {\it amenable in the sense of Bekka} 
\cite{bekka} if there exists a 
state, $\phi$, on the algebra ${\mathcal B}(\H)$
of bounded operators that is $\Ad\, G$-invariant: 
$\phi(\pi(g)T\pi(g)^{-1})=\phi(T)$ for
every $T\in {\mathcal{B}}(\H)$ and every $g\in G$. 
For instance, the homogeneous space $G/H$ is Eymard--Greenleaf amenable if and only if
the quasi-regular representation $\lambda_{G/H}$ of $G$ in $L^2(G/H)$ is 
amenable. 

Let $F$ and $H$ be
closed subgroups of a locally compact group $G$, such that
$F\subseteq H \subseteq G$.
In 1972 Eymard had asked (\cite{Ey}, p. 55):
\begin{itemize}
\item[{\bf Q 1.}]
{\it Suppose the space $G/F$ is amenable. Is then $H/F$ amenable?}
\end{itemize}
This is of course a classical result in the case $F=\{e\}$.

Let $\pi$ be a strongly continuous unitary representation of $H$. In 1990
Bekka has shown \cite{bekka} that, if the unitarily induced representation 
$\ind_H^G(\pi)$ is amenable, then $G/H$ is amenable.
He asked the following:
\begin{itemize}
\item[{\bf Q 2.}]  {\it Does amenability of
$\ind_H^G(\pi)$ imply that of $\pi$?}
\end{itemize}
(This is a more general version of Eymard's question Q 1.)
\begin{itemize}
\item[{\bf Q 3.}] {\it Suppose both $G/H$ and $\pi$ are amenable.
Is then $\ind_H^G(\pi)$ amenable?}
\end{itemize}
Question 3 was partially answered in the positive by
Bekka himself \cite{bekka} under extra assumptions on $H$.
Further discussion can be found in \cite{joli}.  

We show that in general the answer to all three questions above is negative.

\subsection*{\bf Remark} In a recent independent work Nicolas Monod and
Sorin Popa have also solved Eymard's problem (and thus Bekka's Q 2 as well). 
(See their note \cite{MP}
in the same issue.)

\section{Reminders and simple facts}

\subsection{}
A unitary representation $\pi$ of a locally compact group $G$
in a Hilbert space $\H$ is said to {\it almost have invariant vectors}
if for every compact $K\subseteq G$ and every $\e>0$ there is a
$\xi\in \H$ satisfying $\norm \xi=1$ and $\norm{\pi_g(\xi) -\xi}<\e$
for all $g\in K$.

\subsection{\label{alminv}} {\it
If a unitary representation $\pi$ of a group $G$
(viewed as discrete)
almost has invariant vectors, then $\pi$ is amenable.}
\smallskip

This follows from Corollary 5.3 of \cite{bekka}, because
$\pi$ weakly contains the trivial one-dimensional representation, $1_G<\pi$.
Alternatively, choose for every finite $K\subseteq G$ and 
every $\e>0$ an almost
invariant vector $\xi_{K,\e}$ as above, and set
$\phi_{K,\e}(T):=\langle T\xi_{K,\e},\xi_{K,\e}\rangle$. 
Every weak$^\ast$ cluster point, $\phi$, of the net of states $(\phi_{K,\e})$ 
is a $G$-invariant state. 

Note that the converse of \ref{alminv} 
is not true, as for example every unitary 
representation of an amenable group is amenable by Theorem 
2.2 of \cite{bekka}.

\subsection{\label{criterion}} {\it Let $H$ be a closed subgroup of a locally 
compact group $G$. The following statements are equivalent.
\begin{enumerate}
\item[(i)] The homogeneous space $G/H$ is amenable.
\item[(ii)] The left quasi-regular
representation $\lambda_{G/H}$ of $G$ in $L^2(G/H)$, formed with respect
to any (equivalently: every) quasi-invariant measure $\nu$ on $G/H$, almost has
invariant vectors. 
\item[(iii)] The left quasi-regular representation $\lambda_{G/H}$ is
amenable.
\item[(iv)] If $G$ acts continuously by affine transformations on a
convex compact set $C$ in such a way that $C$ contains an $H$-fixed point,
then $C$ contains a $G$-fixed point.
\end{enumerate}
}
The equivalences between (i), (ii), and (iv) are due to 
Eymard \cite{Ey}, while 
(iii) was added by Bekka (Th. 2.3.(i) of \cite{bekka}). 
The condition (ii) is an analogue of Reiter's condition 
$(P_2)$ for homogeneous spaces. (Interestingly, a natural  
analogue of F\o lner's condition for homogeneous spaces fails \cite{Gre}.)

\subsection{\label{extension}}
{\it Let the homogeneous space $G/H$ be amenable, and let $\pi$ be a strongly
continuous unitary representation of $G$.
If the restriction of $\pi$ to $H$ is an amenable representation,
then $\pi$ is amenable as well.}
\smallskip

To prove this statement, denote, following Bekka 
(Section 3 of \cite{bekka}), by $X(\H)$ the
$C^\ast$-subalgebra of ${\mathcal{B}}(\H)$ formed by all operators $T$
with the property that the orbit map of the adjoint action,
\[G\ni g\mapsto \pi(g)T\pi(g)^{-1}\in{\mathcal B}(\H),\]
is norm-continuous. Then $G$ acts upon $X(\H)$ in a strongly 
continuous way by isometries. It follows that the dual adjoint action of $G$
on the state space of $X(\H)$, equipped with the weak$^\ast$ topology,
is continuous. Because of the assumed amenability of $\pi\vert_H$,
there is an $H$-invariant state, $\psi$, on ${\mathcal{B}}(\H)$ and
consequently on $X(\H)$. Eymard's conditional fixed point property
\ref{criterion}.(iv) allows us to conclude
that there is a $G$-invariant state,
$\phi$, on $X(\H)$. This in turn implies amenability of $\pi$
by Theorem 3.5 of \cite{bekka}.

\subsection{} Let $H$ be a closed subgroup of a locally compact group
$G$, and let $\pi$ be a strongly continuous unitary representation
of $H$ in a Hilbert space. By $\ind_H^G(\pi)$ we will denote, as usual, the
unitarily induced representation of $G$. 

\subsection{\label{reciprocity}} The following is
Mackey's generalization of the Frobenius Reciprocity Theorem
(Th eorem 5.3.3.5 of \cite{Wa}). 
\smallskip

{\it Let $H$ be a closed subgroup of a
locally compact group $G$, and let $\pi$ and $\rho$ be finite-dimensional
irreducible unitary representations of $H$ and $G$, respectively.
If $G/H$ carries a finite invariant measure, then 
$\ind_H^G(\pi)$ contains $\rho$ as a discrete direct summand
exactly as many times as $\rho\vert_H$ contains $\pi$ as a discrete direct
summand. 
}
\subsection{\label{stages}}
Let $G$ be a locally compact group, and let $H,F$ be closed subgroups of $G$
such that $F\subseteq H$. Set $\pi=\ind_F^H(1_F)$, where $1_F$ stands
for the trivial one-dimensional representation of $F$. Then
$\pi=\lambda_{H/F}$ is unitarily equivalent to 
the quasi-regular representation of $H$ in $L^2(H/F)$
(\cite{Wa}, Corollary 5.1.3.6).
By the theorem on induction in stages ({\it ibid.}, Proposition 5.1.3.5), 
the induced representation
$\ind_H^G(\lambda_{H/F})$ is unitarily equivalent to $\ind_F^G(1_F)$,
 which is just the representation $\lambda_{G/F}$, 
 the quasi-regular representation of
$G$ in $L^2(G/F)$. 

Now assume that the homogeneous factor space $G/F$ is
amenable. By Bekka's result mentioned above (2.3.(iii)), 
this amounts to the amenability of $\lambda_{G/F}$.

This argument shows that Eymard's Q 1 is a particular case of Bekka's Q 2.

\section{First example}

\subsection{}
Let $H=F_\infty$, where $F_\infty$ denotes the free non-abelian 
group on infinitely many free generators $x_i$, $i\in\Z$.
For every $n\in\Z$, let $\Gamma_n$ denote the normal
subgroup of $F_\infty$ generated by $x_i$, $i\leq n$. 
We will set $F=\Gamma_0$. 

As the factor space $H/F = F_\infty/\Gamma_0$
is a free group, it is not amenable. 

\subsection{} Let the group $\Z$ act on $F_\infty$ by group
automorphisms via shifting the generators:
\[\tau_n(x_m):= x_{m+n},\]
where $m,n\in\Z$ and $\tau$ denotes the action of $\Z$ on $F_\infty$.

Denote by $G=\Z\ltimes_\tau F_\infty$ the semi-direct product formed with
respect to the action $\tau$. That is, $G$ is the Cartesian product
$\Z\times F_\infty$ equipped with the group operation
$(m,x)(n,y)=(m+n,x\tau_m y)$, the neutral element $(0,e)$ and the inverse
$(n,x)^{-1}=(-n,\tau_{-n}x^{-1})$.

Let us show that the homogeneous space $G/F$ is Eymard--Greenleaf amenable.

\subsection{}
For every $n\in\Z$,
\begin{eqnarray*}
(n,e)F&=& \{(n,e)(0,y)\colon y\in F\} \\
&=& \{(n,\tau_n y)\colon  y\in \Gamma_0\} \\
%&=& \{(n,z)\colon z\in \Gamma_n\}\\
&=& \{n\}\times \Gamma_n.
\end{eqnarray*}

Let $S\subset F_\infty$ be an arbitrary 
finite subset. For $n\in\Z$ sufficiently
large, $S\subset \Gamma_n$.

Now for each $(0,s)\in S$,
\begin{eqnarray*}
(0,s)(n,e)F &=& \{(0,s)(n,z)\colon z\in \Gamma_n\} \\
&=& \{n\}\times (s\Gamma_n) \\
&=& \{n\}\times \Gamma_n \\
&=& (n,e)F;
\end{eqnarray*}
that is, the left $F$-coset $(n,e)F\in G/F$ is $S$-invariant. 

Consequently, 
the unit vector $\delta_{(n,e)F}\in \ell^2(G/F)$ is $S$-invariant.
We have proved that the $H$-module $\ell^2(G/F)$ almost has
invariant vectors, and therefore the restriction to $H$ of the left
regular representation $\lambda_{G/F}$ is amenable. 
Since $G/H$ is an amenable homogeneous space ($H$ is normal in $G$ 
and the factor-group $G/H$ is isomorphic to 
$\Z$), we conclude by \ref{extension}
that the left regular representation
of $G$ in $\ell^2(G/F)$ is amenable, that is, that $G/F$ is an amenable 
homogeneous space.

\subsection{} The constructed triple of groups 
$F\subset H\subset G$ provides a negative
answer to the question of Eymard (Q 1). 

In
view of the remarks in \ref{stages}, it provides a negative answer to
Bekka's question (Q 2) as well. Take as $\pi=\ind_F^H(1_F)=\lambda_{H/F}$ the
left quasi-regular representation of $H$ in $\ell^2(H/F)$.
While $\pi$ is not amenable, the induced
representation $\ind_H^G(\pi)=\lambda_{G/F}$ is amenable.

\subsection{} The following way of viewing our example may be instructive. 

If considered as a unitary $F_\infty$-module, $\ell^2(G/F)$
decomposes, up to unitary equivalence, into the orthogonal sum 
\[\ell^2(G/F)\cong \bigoplus_{n\in\Z} \ell^2(F_\infty/\Gamma_n).\]

None of the unitary $F_\infty$-modules $\ell^2(F_\infty/\Gamma_n)$ 
is amenable, yet their orthogonal sum is amenable, because the family 
$(\ell^2(F_\infty/\Gamma_n))_{n\in\Z}$
``asymptotically has invariant vectors'': every finite $S\subset F_\infty$
acts trivially on the Hilbert space $\ell^2(F_\infty/\Gamma_n)$,
provided $n$ is large enough so that $S\subset\Gamma_n$.
(Indeed, $\Gamma_n$ are normal in $F_\infty$.)

\section{Second example}

\subsection{}
Let $G={\mathrm{SL}}(n,\R)$
and $H={\mathrm{SL}}(n,\Z)$, $n\geq 3$. 
Since $H$ is a (non-uniform) lattice
in $G$ ({\it cf.} Exercise 7, \S 2, Chapitre VII of \cite{Bou}), 
the homogeneous space $G/H$ is amenable.

\subsection{}
The group ${\mathrm{SL}}(n,\Z)$ is maximally almost periodic (for instance,
homomorphisms to the groups
${\mathrm{SL}}(n,\Z_p)$, where $p$ is a prime number, separate points
in ${\mathrm{SL}}(n,\Z)$).
Let $\pi$ be a non-trivial (of dimension $>1$)
irreducible unitary finite-dimensional
representation of $H$. Being finite-dimensional, $\pi$ is
an amenable representation; {\it cf.} \cite{bekka}, Theorem 1.3.(i).

\subsection{}
Applying Mackey's Reciprocity Theorem (quoted above, \ref{reciprocity})
with $\rho$ equal to $1_G$, a trivial one-dimensional representation of $G$,
we conclude that the unitarily induced representation
$\ind_H^G(\pi)$ has no non-zero invariant vectors.

At the same time, the Lie group $G={\mathrm{SL}}(n,\R)$, $n\geq 3$,
has property (T) (see {\it e.g.} \cite{dlHV}, 
Chapitre 2.a, Th\'eor\`eme 4), and
the only amenable representations of 
non-compact simple Lie groups with Kazhdan's property $(T)$,
such as $G$, are those
having a non-zero invariant vector. (\cite{bekka}, Remark 5.10;
{\it cf.} also the principal result of \cite{BV}.)

\subsection{}
We conclude: $\ind_H^G(\pi)$ is non-amenable, even if the unitary
representation
$\pi$ and the homogeneous space $G/H$ are both amenable.

This answers in the negative another question posed by Bekka
(Q 3).

\end{document}